\documentclass[11pt]{amsart}
\usepackage{graphicx}
\usepackage{amsmath,amssymb}
\newtheorem{thm}{Theorem}[section]
\newtheorem{cor}[thm]{Corollary}
\newtheorem{lem}[thm]{Lemma}

\theoremstyle{definition}

\theoremstyle{remark}
\newtheorem{rem}[thm]{Remark}
\theoremstyle{Example}
\newtheorem{exa}[thm]{Example}
\numberwithin{equation}{section}

\begin{document}
\title{Congruent numbers via the Pell equation and its analogous counterpart}
\author[Farzali Izadi]{Farzali Izadi }
\address{ Mathematics Department Azerbaijan university of  Tarbiat Moallem ,
 Tabriz, Iran farzali.izadi@gmail.com}
\address{Visitor- Mathematics Department University of Toronto, f.izadi@utoronto.ca}
\begin{abstract}
The aim of this expository article is twofold. The first is to introduce
several polynomials of one variable as well as two variables defined on the positive integers with values as congruent numbers.
The second is to present connections between Pythagorean triples and the Pell equation $x^2-dy^2=1$ plus its analogous counterpart
$x^2-dy^2=-1$ which give rise to congruent numbers n with arbitrarily many prime factors.
\end{abstract}
\maketitle
{\small{\bf Keywords:} Congruent numbers; The Pell equations; Pythagorean triples;\\
 Diophantine equations; elliptic curves}

\section{Introduction}
A right triangle is called rational when its legs and hypotenuse are all rational numbers. Pythagorean triple like $(3,4,5)$ is the most common example of this type of triangles. Scaling such triangle by any rational number $r$, we can get another rational triangle  such as $(3r,4r,5r)$. Any rational right triangle has a rational area but not conversely. A famous longstanding problem about right triangle with rational sides lengths asks: which rational numbers occur as the area of  rational right triangles, i.e., given any rational number $r$, does there exist positive rational numbers $a,b,c$, such that $a^2+b^2=c^2$ and $(\frac{1}{2})ab=r$. This deceptively simple-looking problem is still not completely solved.
In the case of affirmative answer, the number $r$ is called a congruent number, (otherwise it is called non-congruent number ) and a description of all congruent numbers is referred to as the congruent number problem. Since scaling a triangle changes its area by a square factor, and every rational number can be multiplied by a suitable rational square to become a square-free positive integer, we can concentrate in the congruent number problem on square-free positive integers. An equivalent condition for a positive integer $n$ to be a congruent number, as one easily checks, is that there should exist a rational square which, when increased or decreased by $n$, remains a square.

\section{History} The congruent number problem was first stated by the persian mathematician Al-Karaji (C. 953 - C. 1029). However his version involves not in terms of the triangles but instead on its equivalent forms of the square numbers mentioned above. The name comes from the fact that there are three squares which are congruent modulo a number. A major influence on Al-Karaji was the arabic translation of the works of the Greek mathematician Diophantus (C. 210 - C. 290) who posed similar problems.
First tables of congruent numbers were found also in Arab manuscripts namely at the time of Al-Karaji himself, and 5 and 6 appeared there. A small amount of progress was made in the first millennium. In 1225, Fibonacci was the first after Al-Karaji to show that 5 and 7 are congruent numbers. He also stated without giving any proof that 1 is not a congruent number. In years between 1640-1659 Fermat proved this case. The proof was given by the method of descent, which was discovered by Fermat himself on this very problem. For the argument see \cite{bur}. In  1915, Bastein \cite{bas}  found all square-free congruent numbers less than 100 which were exactly 36 cases. Then 62 more square-free congruent numbers less than 1000 were discovered by G\"erardin \cite{ger}.  However by the 1980, there were still cases smaller than 1000 that had not been resolved. The first investigations  on the history of congruent numbers can be found in \cite{lag1}, \cite{lag2}, \cite{guy}, and \cite{dic}. Also, for more information on some other results about congruent and non-congruent numbers without using the theory of elliptic curves consulate \cite{zim}.

\section{Modern results}
 {\it Using elliptic curves:} All recent results about the congruent numbers stem from the fact that n is a congruent number if and
only if the elliptic curve $E_n(Q): y^2 = x^3-n^2x$ contains a rational point with $y\neq 0$, equivalently, a rational
point of infinite order \cite{kob}.
In 1929, Nagell \cite{nag} had a very short and elementary proof of the fact that the rank of $E_n(Q)$ is zero in the case of  $n=p \equiv 3 (mod 8)$
for a prime number $p.$ Thus these numbers are non-congruent. He also pointed out that the same technique shows that 1, 2, and all $2q$ with $q$ a prime $\equiv 5 (mod 8)$ are non-congruent numbers.
Kurt Heegner, exploit the elliptic curves in a nontrivial way by which he  developed the theory of what
is now called Heegner points. In 1952, he proved that all $n=2p$ with prime numbers  $p\equiv 3 (mod 4)$
are congruent. These results are unconditional, while (as described below) most of the
later results rely on the Birch and Swinnerton-Dyer Conjecture \cite{coa}.
Monsky \cite{mon} later extended Heegner method to show that the primes$ \equiv 5, 7(mod 8)$ are congruent numbers. Since primes$ \equiv 3 (mod 8)$ are non-congruent by Nagell's result mentioned above, this only leaves the primes $\equiv 1 (mod 8)$. Here the situation is still unknown. For instance 17 is known to be non-congruent and 41 congruent. In general from the numerical evidence for 1000 numbers in 1972, it has been conjectured by Alter, Curtz, and Kubota \cite{alt1} and \cite{alt2} that if $n \equiv 5, 6, 7 (mod 8)$, then n is a congruent number. It is shown by Stephens \cite{ste}  that this conjecture is a corollary of the Selmer conjecture for elliptic curves . In some special cases, the Selmer conjecture is known to be true by the method of Heegner mentioned above, and  therefore, for these cases ACK conjecture  is also true. These
are when $n$ is a prime congruent to $5$ or $7 (mod 8)$, or when $n = 2p$ where $p$ is a
prime congruent to $3 (mod 4)$.
In 1982, Jerrold Tunnell \cite{tun} made significant progress by exploiting the more trivial connection between congruent numbers and elliptic curves. He found a simple formula for  determining whether or not a number is a congruent number. This allowed the first several thousand cases to be resolved very quickly. One issue is that the complete validity of his formula depends on the truth of a particular case of one of the outstanding problems in mathematics known as the Birch and Swinnerton-Dyer conjecture \cite{coa}.
Let $p,q$ and $r$ denote distinct prime numbers. The following results on congruent numbers have been determined so far.
For congruent numbers:

\begin{itemize}
\item $n=2p_3$ i.e., $p_3\equiv 3(mod 4)$ [Heegner (1952) and Birch (1968)].\
\item $n=p_5,\,p_7$ i.e., $p_i\equiv i(mod 8)$ [Stephen (1975)].\
\item $n=p^uq^v \equiv 5, 6, 7 (mod 8)$, $0 \leq u,v \leq 1$ [B. Groos (1985)].\
\item $n=2p_3p_5, 2p_5p_7$,
\item $n=2p_1p_7$ , $(\frac{p_1}{p_7})=-1$ [P. Monsky (1990)].
\item $n=2p_1p_3$ , $(\frac{p_1}{p_3})=-1$,
\end{itemize}
where $(\frac{p}{q})$ is the Legendre symbol.
A detailed information on non-congruent numbers can be found in \cite{fen}. Also for a recent survey article about the congruent numbers and their variants see \cite{nor}.

In 1984,  Chahal \cite{cha}  applied an identify of Desboves to show that there are infinitely many congruent numbers in each residue class modulo 8 and, in particular, infinitely many square-free congruent numbers, congruent to 1, 2, 3, 4, 5, 6, and 7 modulo 8. This result is generalized later by M. A. Bennet \cite{ben} as follow: If $m$ is a positive integer and $a$ is any integer, then there exist infinitely many (not necessarily square-free) congruent numbers $n$ with $n\equiv a  (mod m)$. If, further, $gcd(a,m)$ is square-free, then there exist infinitely many square-free congruent numbers $n$ with $n\equiv a (mod m)$.

Finally on September 22, 2009 Mathematicians from North America, Europe, Australia, and South America have resolved the first one trillion cases (still conditional on the Birch and Swinnerton-Dyer conjecture) \cite{tri}.

\section{Our results} In the present paper, the author gives an elementary and short proof to construct many congruent 
numbers. First of all, we discuss a very easily verified lemma to show that given any Pythagorean triple, one can find 5 more Pythagorean triples which give rise to some congruent numbers in a very simple manner. Secondly, by using these triples along with the expressions for the initial Pythagorean triple in terms of the two parameters with opposite parity, we obtain 6 different polynomial expressions of two variables defined on the positive integers having values as congruent numbers. Thirdly, according to some simple facts from elementary number theory such as the necessary and sufficient conditions for a number to be written as a sums or differences of two squares, we obtain several other polynomials of one variable which always having values as congruent numbers. Some of these expressions have already been appeared in literature with or without proofs. Finally, in the last part, we will construct simple connections between pythagorean triples and the Pell equation plus its analogous counterpart which give rise to new congruent numbers n with arbitrarily many prime factors.
 To this end, we start off with the following lemma.

\begin{lem}
For a Pythagorean triple $(a,b,c),$ with $b>a$ the following are all Pythagorean triples:

\noindent (1) $(2ac, b^2, a^2+c^2)$, \\
(2) $(2bc, a^2, b^2+c^2)$, \\
(3) $(2ab, b^2-a^2, c^2)$. \
\end{lem}

\noindent{\bf Proof.} It is a straightforward calculation.

\begin{cor}
The numbers $ac$, $bc$, $b^2-a^2$, $a^2+c^2$, and $b^2+c^2$, are all congruent numbers.
\end{cor}
\begin{exa}
For the Pythagorean triple $(3, 4, 5)$ with corresponding congruent number 6 we get 15, 5, 7, 41, 34 as square-free congruent numbers and
20 as congruent number having square factor 4.
\end{exa}

\noindent As is well-known, the primitive Pythagorean triples were completely described by the exact patterns of the form $a=s^2-t^2, b=2st, c=s^2+t^2$, where $s>t\geq 1$ are such that $gcd(s,t)=1,$ and $s,t$ have opposite parity i.e., $s-t$ is always odd.

\begin{cor}
For any positive integer valued parameters $s$ and $t$ with the above conditions, the following expressions  are all congruent numbers.

\noindent (1)$A=st(s^2-t^2)$,\\
(2)$B=st(s^2+t^2)/2$,\\
(3)$C=s^4-t^4$,\\
(4)$D=2(s^4+t^4)$,\\
(5)$E=s^4+t^4+6s^2t^2$,\\
(6)$F=(4s^2t^2-s^4-t^4)$.\\
\end{cor}

\begin{rem}
Another elementary proof for (2) is given in \cite{wad}. By letting $s=x^2$ and $t=2y^2$ in (2) and ignoring square factor, we get
$x^4+4y^2$. This last expression along with (3), (4), (5), and (6) appeared in \cite{alt2} without any explanation or proofs.
\end{rem}
{\bf Example.}
By letting  $s=653821282242$, and $t=127050186481$ in the expression (2), we get the congruent number $n=157*(53156661805)^2.$
Ignoring the square factor we are lead to the square-free congruent number 157.

\noindent The corresponding right triangle is then given by:\\

\noindent

$a=338402045517054238391582296989254448074677078418$,\\

$b=294808091174913744183357386456082152630805118800$,\\

and\

\noindent

$c= 448807035408673939849115026181349726321896944082$.\\

\noindent D. Zagier showed that the  simplest Rational
triangle with area 157 has the following sides.\

\noindent

 $a=6803298487826435051217540/411340519227716149383203$,\\

$b=411340519227716149383203/21666555693714761309610$,\\

and\

\noindent $c=2244035177043369699245755130906674863160948472041$\\

\noindent  divided by\

$8912332268928859588025535178967163570016480830$.\

\noindent Next, we recall  the following elementary lemma.
\begin{lem}
The positive integer $n$ can be written as a difference of two squares namely
$n=a^2-b^2$ if and only if $n-2$ is not divisible by 4.
\end{lem}
For the proof of this  lemma the reader may consulate \cite{ros}.
This simple fact along with the above corollaries lead to the following results.

\begin{thm}
For any positive integer $k$, the following are all congruent numbers.
\begin{eqnarray*}
&&(A)\left\{
     \begin{array}{ll}
       (A1)\,\,k(k^2-1), & \hbox{} \\
       (A2)\,\,2k(k^2+1), & \hbox{$for \,\,\, n=4k;$} \\
       (A3)\,\,(k^4-1), & \hbox{}
     \end{array}
   \right.\\
&&(B)\left\{
     \begin{array}{ll}
       (B1)\,\,k(2k+1)(4k+1), & \hbox{} \\
       (B2)\,\,2k(2k+1)(8k^2+4k+1), & \hbox{$for\,\,\, n=4k+1;$} \\
       (B3)\,\,(4k+1)(8k^2+4k+1), & \hbox{}
     \end{array}
   \right.\\
&&(C)\left\{
  \begin{array}{ll}
    (C1)\,\,k(k+1)(2k+1), & \hbox{} \\
    (C2)\,\,2k(k+1)(2k^2+2k+1), & \hbox{$for\,\,\, n=4k+2 \equiv n=2k+1;$} \\
    (C3)\,\,(2k+1)(2k^2+2k+1), & \hbox{}
  \end{array}
\right.\\
&&(D)\left\{
     \begin{array}{ll}
       (D1)\,\,2(k+1)(2k+1)(4k+3), & \hbox{} \\
       (D2)\,\,(k+1)(2k+1)(8k^2+12k+5), & \hbox{$for\,\,\, n=4+3;$} \\
       (D3)\,\,(4k+3)(8k^2+12k+5), & \hbox{}
     \end{array}
   \right.\\
&&(E)\left\{
     \begin{array}{ll}
       (E1)\,\,2k(2k+1)(2k-1), & \hbox{} \\
       (E2)\,\,k(4k^2+1), & \hbox{$for\,\,\, n=8k;$} \\
       (E3)\,\,(2k+1)(4k^2+1), & \hbox{}
     \end{array}
   \right.\\
\end{eqnarray*}
\begin{eqnarray*}
&&(F)\left\{
     \begin{array}{ll}
       (F1)\,\,k(4k+1)(8k+1), & \hbox{} \\
       (F2\,\,2k(4k+1)(32k^2+8k+1), & \hbox{$for\,\,\, n=8k+1;$} \\
       (F3)\,\,(8k+1)(32k^2+8k+1), & \hbox{}
     \end{array}
   \right.\\
&&(G)\left\{
     \begin{array}{ll}
       (G1)\,\,2(2k+1)(4k+1)(8k+3), & \hbox{} \\
       (G2)\,\,(8k+3)(32k^2+24k+5), & \hbox{$for\,\,\, n=8k+3;$} \\
       (G3)\,\,(2k+1)(4k+1)(32k^2+24k+5), & \hbox{}
     \end{array}
   \right.
\end{eqnarray*}
\end{thm}

\noindent {\bf Proof.} For the proof, we assume that $n$ is an odd number. For an even number, the same arguments can be easily
applied. Let $n=2mk+(2r+1)$, where $m=1$ or $m$ is an even number, and $s$ varies from $0$ to $m-1.$ Then we can write $n$ as\\

$n=2mk+(2r+1)=(mk+r+1)^2-(mk+r)^2.$\\

\noindent  From this we define a Pythagorean triple as follows:

\noindent $a=s^2-t^2$, $b=2st$, and $c=s^2+t^2$, where $s=(mk+r+1)$, and $t=(mk+r)$. Now, the results are immediate from the above corollary.

\begin{rem}
From (A1), we see that the product of any three consecutive numbers is a congruent number.
 In particular, the product of twin primes are parts of congruent
numbers. If moreover, one of the three numbers appearing in the product is square, then
the product of the two remaining numbers is a congruent number.
\end{rem}

\begin{rem}
(E2) has already been discussed by Chahal \cite{cha}. In his paper, he obtained this result as a consequence of
 some properties of elliptic curves along with an identity of Desboves.
\end{rem}

Before stating the main theorems, we recall the following well-known fact from elementary number theory.
The proof can be found in \cite{ros}.
\begin{lem}
The positive integer $n$ is the sum of two squares if and only if each prime factor $n$ of the form
$4r+3$ occurs to an even power in the prime factorization of $n$.
\end{lem}

\begin{thm}
Let $d=p_{1}p_{2}...p_{m}$, where all $p_{i}s$ are distinct primes such that $p_{1}=2$, and all other
primes are of the form $4r+1$. Then for any $x$-component of the solution
 $(x,y)$  in the Pell-like equation $x^2-dy^2=-1$,
 $2xd$ is a congruent number.
Furthermore, if $x$ a square itself, then $2d$ is also a congruent number.
\end{thm}
\noindent {\bf Proof.} For any $x$-component of the solution , the triple $(x^2-1, 2x, x^2+1)$ is a Pythagorean. Now main lemma implies
that $2x(x^2+1)$ is a congruent number. Substituting $x^2+1$ by $dy^2$ from the Pell-like equation and ignoring the square factor $y^2$, we see that
$2xd$ is a congruent number. The last assertion is a triviality.
\begin{cor}
For any number $d$ as in the statement of the theorem, there exists a congruent number having $d$ as a square-free part with arbitrarily many
prime factors of the above mentioned forms.
\end{cor}

\begin{thm}
Let $d$ be any positive number which is not in the form $4r+2$. Then for any $x$-component of the solution in the pell equation
$x^2-dy^2=1$ other than $(\pm{1},0),$
$xd$ is a congruent number. Furthermore, if $x$ is a square itself, then $d$ is also a congruent number.
\end{thm}

\noindent {Proof.} In this case, we consider the number $x(x^2-1)$ which is congruent by the same reasoning. Replacing $x^2-1$ by $dy^2$, and
ignoring the square factor $y^2$, the result follows easily.

\begin{rem}
It is well known that for any positive square-free integer d, the Pell equation $x^2-dy^2=1$ has an
infinitude of solutions, which can be easily expressed in terms of the fundamental
solution p, q, where $p, q > 0$. On the other hand, its analogous counterpart with equation $x^2-dy^2 = -1$
is solvable for only certain values of $d$. In fact, for this equation a necessary condition to be solvable is that all odd prime factors of $d$ must be of the form  $4r+1$, and that  cannot be doubly even (i.e., divisible by 4). However, these conditions are not sufficient for a solution to exist, as demonstrated by the equation $x^2-34y^2=-1$ , which has no solutions in integers, to see this the reader may consulate \cite{nag2}.
\end{rem}
{\bf Example.} Let $d=m^2+1$. Then the Pell equation $u^2-dv^2=1$ has a solution $(2m^2+1, 2m)$. From this it is easy to see that for any nonzero $y$-component in the solution of the Pell equation $x^2-2y^2=1$, the number $d=y^2+1$ is a congruent number.


\begin{thebibliography}{99}
\bibitem{alt1} Alter, R., Curtz, T. B., Kubota, K. K., : Remarks and results on congruent numbers.
Proc. Third Southeastern Conf. on Combinatorics, Graph Theory and Computing 1972, pp. 27--35.

\bibitem{alt2} Alter, R., Curtz, T. B., : A note on congruent numbers. Math. Comp. <(1974), 303--305.
\bibitem{tri} http://www.aimath.org/news/congruentnumbers/
\bibitem{bas} Bastein, L., :Nombers congruents. Interm\"ediaire des Math. 22 1915, 231-232.
\bibitem{ben} Bennett, M. A., :Lucas' square pyramid problem revisited, Acta Arith. 105 (2002) 341-347.
\bibitem{bur} Burton, D. M., :Elementary number theory, 6th ed., McGaw-Hill, New York, 2007.
\bibitem{cha} Chahal, j. S., :On an identity of Desboves Proc. Japan Acad. Ser. A Math. Sci. V.60, No. 3 (1984), 105-108.
\bibitem{coa} Coates, J., Wiles, A., : On the conjecture of Birch and Swinnerton-Dyer. Invent. math. (1977), 223--251.
\bibitem{dic} Dickson, L. E., :book History of the Theory of Numbers Volume II (ISBN 0-8218-1935-6) Chapter XVI.
\bibitem{fen} Feng, K., :Non-congruent numbers, odd graphs and the Birch—Swinnerton-Dyer conjecture, Acta Arith. 75 (1996), 71–83.
\bibitem{ger} G\"erardin, A., :Nombers congruent. Interm\"ediaire des Math. 22 1915, 52-53.
\bibitem{guy} Guy, R., :Unsolved Problems in Number Theory, 2nd ed. Sprniger Verlag 1994.
\bibitem{kob} N. Koblitz: Introduction to Elliptic Curves and Modular Forms, 2nd ed. New York: Springer-Verlag 1993.
\bibitem{lag1} Lagrange, J., : Nombres congruents et courbes elliptiques. Sém. Delange-Pisot-Poitou (Théorie des nombres) 16e année, 1974/75,
 no. 16, 17 p.
\bibitem{lag2} Lagrange, J., : Construction d'une table de nombres congruents. Manuscript, Reims.
\bibitem{htt} http://mathworld.wolfram.com/CongruentNumber.html
\bibitem{mon} Monsky, P., :Mock Heegner points and congruent numbers, Math. Z. 204 (1990),  45–68.
\bibitem{nag}T. Nagell, L' Analyse ind´etermin´ee de degr´e sup´erieur, vol. 39, Gauthier-Villars, Paris, 1929  16-17.

\bibitem{nag2} Nagell, T., :Introduction to Number Theory. New York: Wiley, pp. 195-212, 1951. §56-58, "The Diophantine Equation $x^2-dy^2=\pm{1}$ ."
\bibitem{wad} K. Noda and H. Wada, All congruent numbers less than 1000, Proc. Japan Acad. 69 (1993), p. 175–178.
\bibitem{ros}Rosen, K., :Elementar number theory and its applications, 5th ed. Eddison-Wesley 2005.
\bibitem{ste}N. M. Stephens: Congruence properties of congruent numbers, Bull. London Math. Soc. 7 (1975), 182–184.
\bibitem{nor} Top, J.,  Noriko, Y. :Congruent numbers and their invariants, Algorithmic NT, MSRI Pub. V.44, 2008.
\bibitem{tun} Tunnell, Jerrald B. A :classical Diophantine problem and modular forms of weight 3/2. Inv. Math. 72 (2): 323–334.
\bibitem{zim} Zimmer, H. G., Congruent numbers-From elementary to alg. NT ACM SIGSAM Bulletin archive
V. 18-19, Issue 4-1  (November-February 1984-1985)  54-62.
\end{thebibliography}
\end{document}